\documentclass[12pt,leqno]{amsart}
\usepackage{tikz}
\usetikzlibrary{positioning}

\makeatletter
\def\l@subsection{\@tocline{2}{0pt}{2pc}{6pc}{}}
\makeatother

\usepackage{comment}

\usepackage{amsmath}
\usepackage{amssymb}
\usepackage{mathtools}
\usepackage{amsmath}
\usepackage{amssymb}
\usepackage{latexsym}
\usepackage{amsopn}
\usepackage{amsthm}
\usepackage{url}
\usepackage{multicol}
\usepackage[
  paperwidth=6in,
  paperheight=9in,
  textwidth=4.7in,
  textheight=7.5in,
  centering
]{geometry}

\usepackage{graphicx}
\graphicspath{ {./images/}}

\usepackage{tikz}
\usetikzlibrary{positioning}

\usepackage{accents}

\def\underTilde#1{{\baselineskip=0pt\vtop{\hbox{$#1$}\hbox{$\sim$}}}{}}

\def\bSigma{\underTilde{\Sigma}}
\def\bPi{\underTilde{\Pi}}
\def\bDelta{\underTilde{\Delta}}

 \newtheorem{thm}{Theorem}[section]
 
 \newtheorem{lem}[thm]{Lemma}

\newtheorem*{theorem1}{Theorem 1}
\newtheorem*{theorem2}{Theorem 2}
\newtheorem*{theorem3}{Theorem 3}
\newtheorem*{theorem4}{Theorem 4}

\newtheorem{prop}[thm]{Proposition}

\theoremstyle{definition}
\newtheorem{defn}[thm]{Definition}

\newtheorem{rem}[thm]{Remark}

\theoremstyle{remark}

 \numberwithin{equation}{section}

\makeatletter
\renewcommand\thepart{\Roman{part}}

\renewcommand\part[1]{%
  \cleardoublepage
  \refstepcounter{part}%
  \addcontentsline{toc}{part}{Part \thepart: #1}%
  \begin{center}
    \bfseries\Large Part \thepart\\[1ex]
    \bfseries #1
  \end{center}
  \vspace{2em}
}
\makeatother

\begin{document}


\title{Universally Measurable Sets and Correlated Equilibria}

\author{Stuart Zoble \\ November 9, 2026}

\thanks{{\em Acknowledgments.}  We thank Adam Brandenburger for many illuminating discussions about game theory over the years.}

\address{Department of Mathematics, Albert Nerken School of Engineering, The Cooper Union for the Advancement of Science and Art, New York City}
\email{szoble@gmail.com \\ stuart.zoble@cooper.edu}


\begin{abstract} \sloppy We study alternate definitions of correlated equilibrium in games with Polish strategy spaces.  We first consider universally measurable decision problems equipped with a measure assignment and a prior, comparing deviation stability with the stronger requirement of conditional rationality.  We show they are equivalent for Borel measurable problems, but their equivalence for all analytically measurable problems depends on the model of set theory.  The full equivalence follows from $\mathsf{MUP}$, while a nonmeasurable set of reals yields a counterexample.  When the set of strategies which are best responses to their assigned measures is measurable by the prior, we characterize the equivalence by the principle $\mathsf{MUP}^{\mathsf{UM}}$, which holds in a model of Larson--Neeman--Shelah.  The results carry over to games, where we further show that conditionally rational priors must concentrate on $\text{RCBR}$ strategies, while deviation stable priors may fail to do so even in games with analytically measurable payoffs.  This provides an epistemic justification for defining a correlated equilibrium as a conditionally rational prior. \end{abstract}


\maketitle




\section{Introduction}  A set-theoretic perspective on correlated equilibrium, Aumann's classical generalization of Nash equilibrium, reveals a distinction between two natural formulations that is obscured by the restriction to finite or continuous games.

Consider a Nash equilibrium $(P_1,P_2)$ of a game $G = (\pi_1,\pi_2)$, with each $$\pi_i:X_1 \times X_2 \rightarrow \mathbb{R}$$
a bounded, universally measurable map. Thus each $P_i \in \Delta(X_i)$ and for each $Q_1 \in \Delta(X_1)$ and $Q_2 \in \Delta(X_2)$, $$\pi_1(Q_1,P_2) \leq \pi_1(P_1,P_2) \mbox{ and } \pi_2(P_1,Q_2) \leq \pi_2(P_1,P_2),$$ where we extend each $\pi_i$ by taking expectations.  
To generalize this concept to an arbitrary joint distribution $P \in \Delta(X_1 \times X_2)$, let us try to capture the Nash property in terms of the measure $$P = P_1 \otimes P_2 \in \Delta(X_1 \times X_2).$$  The requirement that  
$$\int_{X_1 \times X_2} \pi_1( s^*, t) dP(s,t) \leq \int_{X_1 \times X_2} \pi_1( s, t) dP(s,t)$$ for any $s^* \in X_i$, and similarly for player $2$, is easily seen to be equivalent to the Nash property when $P = P_1 \otimes P_2$.  For general $P$, we allow the deviation to depend on $s$ so we require $$\int_{X_1 \times X_2} \pi_1( d(s), t) dP(s,t) \leq \int_{X_1 \times X_2} \pi_1( s, t) dP(s,t)$$ for any measurable $d: X_1 \rightarrow X_1$, and similarly for player $2$.
This is the definition of Aumann in \cite{aumann1} and \cite{aumann2}, where the concept of a correlated equilibrium was introduced.  We will call this property \emph{deviation stability}.

Looking at the Nash pair $(P_1,P_2)$ again, we see that $$P_i ( \{ s \in X_i \ | \ s \mbox{ is a best response to } P_j \}) = 1$$ for each $i$, that is each $P_i$ concentrates on best responses to $P_j$ (strategies that yield as high a payoff as any alternative in $X_i$ against $P_j$).  For an arbitrary joint distribution $P$, we would require that $$\mbox{marg}_{X_i}(P)  ( \{ s \in X_i \ | \ s \mbox{ is a best response to } \mu_i^P(s) \}) = 1$$ where $$\mu_i^P:X_i \rightarrow \Delta(X_j)$$ is part of some (equivalently any) disintegration of $P$, that is any pair of Borel measure assignments $(\mu_1^P, \mu_2^P)$ with $$P(B_1 \times B_2) = \int_{B_1} \mu_1^P(s)(B_2) dP_1(s) = \int_{B_2} \mu_2^P(t)(B_1) dP_2(t)$$
defining $P$ for all Borel sets $B_i \subseteq X_i$, where each $P_i = \mbox{marg}_{X_i}(P)$.  We call such a $P$ \emph{conditionally rational}. This definition can be found for example in Yildiz's \emph{Graduate Game Theory} and is easily seen to imply deviation stability.

In our general setting of universally measurable maps on Polish spaces, we find that the two properties coincide for Borel measurable payoffs, but depending on the model of set theory, may or may not coincide for analytically measurable payoffs.  The full equivalence follows from the measure uniformization principle\footnote{$\mathsf{MUP}$ holds in the Solovay model and is a consequence of $\mathsf{AD}$ (see \cite{measuregame}).  We will show that $\mathsf{MUP}^{\mathsf{UM}}$ holds in a $\mathsf{ZFC}$ model of Larson--Neeman--Shelah from \cite{LNSh}.}, $\mathsf{MUP}$, and is exactly characterized by $\mathsf{MUP}^{\mathsf{UM}}$ for games and priors where  $$\{ s \in X_i \ | \ s \mbox{ is a best response to } \mu_i^P(s) \}$$ is measurable by $\mbox{marg}_{X_i}(P)$ for each $i$.  This condition rules out a counterexample created from a bijection between a nonmeasurable set and a universally null set of reals.  We develop these results in the context of decision problems, and they carry over immediately to games.  

\begin{theorem1}  $\mathsf{MUP}$, and hence $\mathsf{AD}$, imply that deviation stability and conditional rationality coincide.  In $\mathsf{ZFC}$ there is a counterexample for universally measurable games.
\end{theorem1}

\begin{theorem2}  For rationality-measurable common strategy priors, the 
equivalence of deviation stability and conditional rationality is equivalent to $\mathsf{MUP}^{\mathsf{UM}}.$
\end{theorem2}

Finally, we show that a correlated equilibrium (in the conditionally rational sense), like the Nash equilibrium it generalizes, must concentrate on strategies that are consistent with rationality and common belief in rationality ($\mathrm{RCBR}$ strategies), whereas a deviation stable strategy prior may not, even in the analytically measurable case. This provides an epistemic justification for the one formulation over the other.

\begin{theorem3}  A common strategy prior is conditionally rational if and only if it is the strategy marginal of a common prior concentrated on $\mathrm{RCBR}$ states.  In particular, every conditionally rational common strategy prior concentrates on $\mathrm{RCBR}$ strategies.\footnote{For the corresponding chracterization in finite games see Dekel and Siniscalchi \cite[Theorem 2.14]{dekel} building on Aumann \cite{aumann2}. }
\end{theorem3}

\begin{theorem4}  It is independent of $\mathsf{ZFC}$ whether every deviation stable common strategy prior in an analytically measurable game concentrates on $\mathrm{RCBR}$ strategies.  \end{theorem4}

Thus, by Thoerem 3, whether or not a deviation stable prior concentrates on RCBR strategies, if it is not conditionally rational, then it cannot be the strategy marginal of any common state space prior concentrated on RCBR states.

\section{Preliminaries} We begin with a few observations about the $\sigma$-algebra of universally measurable sets, $\Gamma^{\mathsf{UM}}$, and a version of Solovay's measure uniformization principle for these sets.   Throughout we will identify Borel probability measures with their completions.  

\begin{defn} A subset $A$ of a Polish space $X$
is universally measurable if it is measurable in the completion of every $\mu \in \Delta(X)$.  $A$ is universally null if it has measure zero in the completion of every atomless $\mu \in \Delta(X)$.  We let $\Gamma^{\mathsf{UM}}(X)$ denote the $\sigma$-algebra of universally measurable subsets of $X$ which is closed under Borel preimages.
\end{defn}

\begin{defn}  We let $\mathsf{MUP}^{\mathsf{UM}}$ denote the principle asserting that any $\Gamma^{\mathsf{UM}}$ relation $R \subseteq 2^{\omega} \times 2^{\omega}$ with $\nu(\mbox{proj}(R))>0$ can be uniformized almost everywhere on $\mbox{proj}(R)$ by a Borel map, where $\nu$ is the standard coin-flipping measure on $2^{\omega}$. \end{defn}

\begin{rem} We obtain an equivalent principle by requiring that 
$$\nu(\mbox{proj}(R)) = 1$$ in the definition above, or by requiring that $R$ is total.  This principle implies the analogous principle for any uncountable Polish space and Borel probability measure.  \end{rem}

\begin{prop}[Relevant Results for $\Gamma^{\mathsf{UM}}$] 
\leavevmode
\begin{enumerate}
    \item \textnormal{(Kechris)} $\bPi^1_1$ determinacy implies $\bSigma^1_2 \subseteq \Gamma^{\mathsf{UM}}$
    \item $\textnormal{(Martin, Solovay)}$  $\mathsf{MA} + \neg \mathsf{CH}$ implies $\bSigma^1_2 \subseteq \Gamma^{\mathsf{UM}}$
    \item \textnormal{(Solovay)} $(\bDelta^1_2)^{\text{ZF}} \subseteq \Gamma^{\mathsf{UM}}$
    \item \textnormal{(Larson--Shelah)} It is consistent that $\Gamma^{\mathsf{UM}} \subsetneq \bDelta^1_2$ 
    \item \textnormal{(Grzegorek)} $\Gamma^{\mathsf{UM}}$ is not closed under projections
    \item $\mathsf{MUP}^{\mathsf{UM}}$ implies $\bDelta^1_2 \subseteq \Gamma^{\mathsf{UM}}$. 
    \item $\mathsf{MUP}^{\mathsf{UM}}$ fails in $\mathsf{L}$.  
\end{enumerate}
\end{prop}
\begin{proof} For (1) see 6G.12(a) of Moschovakis,  attributed to Kechris.  To see (6) suppose $A \in \bDelta^1_2$ so that $A = \mbox{proj}(P_0)$ and $X \setminus A = \mbox{proj}(P_1)$ where $P_0,P_1$ are coanalytic.  Define $$R(x,(i,y)) \leftrightarrow P_i(x,y)$$ which is coanalytic and hence $\Gamma^{\mathsf{UM}}$.  A uniformization of $R$ (replaced by a Borel function that agrees almost everywhere with respect to the measure under consideration) gives the desired Borel witness to the measurability of $A$.  Thus $\mathsf{MUP}^{\mathsf{UM}}$ fails in the Larson--Shelah model (a forcing extension of $\mathsf{L}$) by the fact that the containment in (4) is proper (see Remark 6.3 of \cite{LSh}).

In $\mathsf{L}$ there is a coanalytic total function $$f:2^{\omega} \rightarrow 2^{\omega}$$ with a thin graph $R$ (see Moschovakis 5A.6), that is, $R$ does not contain a perfect subset.  To sketch the definition of $f$, we use the fact that there is a $\bPi^1_1$ set $Q(x,y,w)$ such that $\exists w \ Q(x,y,w)$ holds 
if and only if $x <_{\mathsf{L}} y$ and $y$ is the 
$<_{\mathsf{L}}$-least code of some ordinal.  The function $f$ is a coanalytic uniformization which outputs the unique real coding the unique pair $(y,w)$ such that $Q(x,y,w)$ holds.  By assumption $f$ is a total function.  Any closed subset of $$R = \mbox{graph}(f)$$ would have analytic projection on the second coordinate, which would therefore be countable by boundedness.  Thus $R$ is universally measurable and has full projection on the first coordinate.  If $\mathsf{MUP}^{\mathsf{UM}}$ held there would be an Borel selector $B$ and an uncountable Borel set $B$ with $$\mbox{graph}(g \upharpoonright B) \subseteq R,$$ which must contain a perfect set, contradicting thinness, and establishing (7).

To prove (5), that is, $$\neg (\exists^{\mathbb{R}} \Gamma^{\mathsf{UM}} \subseteq \Gamma^{\mathsf{UM}}),$$ we use the 
graph of a function with nonmeasurable domain and universally null range.  Let $\lambda$ be the coin-flip measure on $2^{\omega}$.  Using a measure preserving Borel isomorphism to move to $2^{\omega}$, we can apply a theorem of Grzegorek (see Theorem 2.8 of \cite{LNSh}) to find a universally null set of the same cardinality as a nonmeasurable set. Thus we may choose a nonmeasurable set
$A\subseteq 2^{\omega}$, a universally null set
$N\subseteq 2^{\omega}$, and a bijection $f:A\to N$.
Write
\[
F=\{(s,f(s)) \ | \ s\in A\}.
\]
We first observe that $F$ is universally null.
Indeed, let $\rho$ be an atomless Borel probability
measure on $2^{\omega} \times 2^{\omega}$, let $\nu$ be its second marginal,
and let $D$ be the countable set of atoms of $\nu$.
Since $\nu$ restricted to $2^{\omega} \setminus D$ is atomless,
universal nullity of $N$ gives
\[
\nu(N\setminus D)=0.
\]
Consequently, $$\rho( F \cap ( 2^{\omega} \times (N \setminus D))  ) = 0.$$
The remaining part is countable,
because $f$ is injective, and is therefore also
$\rho$-null. Thus $\rho(F) = 0$, and it follows that 
$F$ is universally null and hence $F \in \Gamma^{\mathsf{UM}}$.  But $A = \mbox{proj}(F) \notin \Gamma^{\mathsf{UM}}$.

We mention (3) as the provably $\bDelta^1_2$ sets are a  large fragment of 
$\bDelta^1_2$, and refer to \cite{fen} for details.

That $\mathsf{MA} + \neg \mathsf{CH}$ implies every $\bSigma^1_2$ set is universally measurable is due to Martin and Solovay \cite{MS}.  To see this, note that every $\bSigma^1_2$ set is $\omega_1$-Borel.  Under $\mathsf{MA} + \neg \mathsf{CH}$, the union of $\omega_1$ $\mu$-null sets is $\mu$-null for any $\mu \in \Delta(2^{\omega})$ (see pages 167-169 of \cite{MS}).  So let $B$ be the union of Borel sets $\{B_{\alpha} \ | \ \alpha < \omega_1\}$.  There is $\beta$ with $$\mu( \bigcup_{\alpha<\beta} B_{\alpha}) = \mbox{sup}_{\gamma < \omega_1} \ \mu(\bigcup_{\alpha<\gamma} B_{\alpha})$$ so
we must have $$\mu( B_{\gamma} \setminus \bigcup_{\alpha<\beta} B_{\alpha}) = 0$$ for $\gamma \geq \beta$, and hence 
$$\mu(B \setminus \bigcup_{\alpha<\beta} B_{\alpha}) = 0,$$ which shows that $B$ is measurable.  This proves (2).  \end{proof}

\begin{prop} $\mathsf{MUP}^{\mathsf{UM}}$ holds in the random model.
\end{prop}

\begin{proof}  We work in the model of Larson-Neeman-Shelah \cite{LNSh}, which we can take to be a forcing extension of $\mathsf{L}$ obtained by adding $\omega_2$ random reals, which can be referred to as the random model.  By Theorem 2.7 of that paper, and the discussion preceding their Theorem 2.10, 
\begin{enumerate}
    \item[(A)] Every $\Gamma^{\mathsf{UM}}$ set is $\omega_1$ Borel, and 
    \item[(B)] $\mbox{cov}(\mathcal{N}) = \omega_2$,
\end{enumerate}
meaning no positive measure Borel set can be covered by $\omega_1$ null sets (see also Theorem 1.1 of \cite{LNSh}; the normalized restriction of the coin-flip measure to such a Borel set is isomorphic to Lebesgue measure on the unit interval).  This puts us in a situation similar to (2) of the last proposition.

Suppose $R \subseteq 2^{\omega} \times 2^{\omega}$ is a relation with full projection which is universally measurable.  Let $\mu$ be the coin flip measure.  By (A) $$R = \bigcup_{\alpha<\omega_1} B_{\alpha}$$ where each $B_{\alpha}$ is Borel and so each $$A_{\alpha} = \mbox{proj}(B_{\alpha}) \in \bSigma^1_1$$ is universally measurable, and their union exhausts $2^{\omega}$.  Let $M$ be the supremum of measures of unions of countable subfamilies of $\{A_{\alpha} \ | \ \alpha < \omega_1\}$, and let $\{ A_{\alpha} \ | \ \alpha < \beta \}$ realize this supremum. Let $$H = \bigcup_{\alpha <\beta} A_{\alpha}.$$  We argue that $\mu(H) = 1$.  For every $\alpha \geq \beta$, the maximality of $M$ gives $\mu(A_{\alpha} \setminus H) = 0$.  If $\mu(H) < 1$ choose a positive measure Borel $C \subseteq 2^{\omega} \setminus H$.
Since $R$ has full projection $$C = \bigcup_{\alpha \geq \beta} C \cap A_{\alpha}$$ which contradicts (B).  Hence $\mu(H) = 1$.  It follows that the Borel relation $$B = \bigcup_{\alpha < \beta} B_{\alpha} \subseteq R$$ has conull projection $H$.  By the Jankov--von Neumann theorem, $B$ admits an analytically measurable uniformizing map defined on $H$, which agrees with a Borel map on a conull set. This Borel map is the desired uniformization of $R$. \end{proof}
.

\section{Decision Problems}

The set-theoretic issues of interest here can be understood in the decision theoretic context.

\begin{defn}[Decision Problems] We call a bounded universally measurable map $$\pi: X \times Y \rightarrow \mathbb{R}$$ on Polish spaces $X$ and $Y$ a \emph{decision problem}.  We add some additional structure to $\pi$, namely a prior $P \in \Delta(X)$ and a Borel measure assignment $$\mu:X \rightarrow \Delta(Y),$$ or belief assignment, which must be defined $P$-almost everywhere. We call $(\pi,\mu,P)$ a \emph{decision problem with beliefs and a prior}. \end{defn}
\begin{rem} We extend such a $\pi$  
to $$\pi:\Delta(X) \times \Delta(Y) \rightarrow \mathbb{R}$$ by taking expectations, and write $\pi(\mu,\nu)$ for $\mu \in \Delta(X)$ and $\nu \in \Delta(Y)$ for such an expectation.  \end{rem}

\begin{defn}[Deviation Stability] The problem $(\pi,\mu,P)$ is \emph{deviation stable} if there is no Borel measurable map $$d:X \rightarrow X$$ such that 
$$\int_{X} \pi(x,\mu(x)) dP(x) < \int_{X} \pi(d(x),\mu(x)) dP(x).$$  \end{defn}

\begin{rem}  An equivalent definition is obtained by allowing deviations to be $P$-measurable, or universally measurable.  \end{rem}

\begin{defn}[Conditional Rationality] We say that the problem is \emph{conditionally rational} if $x$ is a best response to $\mu(x)$ for $P$-almost every $x$, that is for $x$ on a set of $P$-measure one $$\pi(x,\mu(x)) \geq \pi(y,\mu(x))$$ for every $y \in X$. \end{defn}

\begin{prop}  Let $(\pi,\mu,P)$ be a decision problem that is conditionally rational.  Then it is deviation stable.  \end{prop}

\begin{proof} Immediate.  \end{proof}

We are interested in the converse.  Recall that $\pi$ is universally measurable by assumption.  The issue is that a failure of conditional rationality means that the set
$$B = \{ (s,a) \ | \ \pi(a,\mu(s)) > \pi(s,\mu(s))  \}$$ has nonnull projection.  While $B$ is universally measurable, the two issues are 
\begin{enumerate}
    \item[(P1)] the measurability of $\mbox{proj}(B)$, and 
    \item[(P2)] finding a Borel uniformization of $B$,
\end{enumerate}
\noindent 
which are resolved in the three cases below. 
\begin{prop}[From Uniformization]  Let $(\pi,\mu,P)$ be a decision problem that is deviation stable.
\begin{enumerate}
    \item If $\pi$ is Borel then it is conditionally  rational.
    \item Assuming $\bPi^1_1$-determinacy, if $\pi$ is analytically measurable then it is conditionally  rational.
    \item Assuming $\mathsf{MUP}$, $\pi$ is conditionally rational.
\end{enumerate}
\end{prop}
\begin{proof} Let $B$ be the strict improvement relation $$B = \{ (s,a) \ | \ \pi(a,\mu(s)) > \pi(s,\mu(s))  \}.$$ 
For (1), the set $\mbox{proj}(B) \in \bSigma^1_1$ is measurable and so $B$ can be uniformized everywhere by Jankov--von Neumann, and that map can then be replaced by a Borel map agreeing with it almost everywhere. 

For (2), we argue that $B \in \bSigma^1_2$ and so 
$$\mbox{proj}(B) \in \bSigma^1_2,$$ and the Proposition 2.4(1) applies and we can similarly use the Kondo-Novikov-Addison uniformization theorem.  For this complexity calculation, it suffices to show that the extension of $\pi$ to $X \times \Delta(Y)$ is analytically measurable given our assumption that $\pi$ is bounded and analytically measurable.  For analytic $A \subseteq X \times Y$ the map from $(s,\nu)$ to $\nu(A_s)$ is analytically measurable by Theorem 29.26 in Kechris \cite{kechris}.  A monotone class argument extends this conclusion to analytically measurable sets $A$, and the result follows for $\pi$ by a simple function approximation argument.    

For (3), $\mathsf{MUP}$ implies all sets have measurable uniformizations (see \cite{Zoble}).
\end{proof}
\begin{defn}[Decision Problem of a Relation]  Given a universally measurable relation $R \subseteq 2^{\omega} \times 2^{\omega}$ we set $$X =\{0,1\} \times 2^{\omega} \mbox{ and }Y = 2^{\omega}$$ and define the decision problem $\pi^R$ on $X \times Y$ by $$\pi^R((1,y),x) = \chi_R(x,y)$$ and $\pi^R((0,y),x) = 0$.  We will use this problem with the measure assignment $$\mu((0,s)) = \delta_s,$$ and various priors $P$ which concentrate on $\{0\} \times 2^{\omega}$.  For a decision problem of this form, a strategy $s$ is rational if and only if $s \notin \mbox{proj}(R)$, and a profitable Borel deviation gives rise to a partial uniformization of $R$ defined on a set of positive $P$-measure.  \end{defn}

\begin{prop}[Uniformization Necessary] If every deviation stable problem is conditionally rational then 
$\mathsf{MUP}^{\mathsf{UM}}$ holds.
\end{prop}
\begin{proof} Let $R$ be an arbitrary total relation and $\pi^R$ the decision problem from Definition 3.8. 
For any Borel set of positive measure $B \subseteq 2^{\omega}$ let $\nu_{B}$ be the coin-flip measure normalized and regarded as being on $\{0,1\} \times 2^{\omega}$ with $\nu_B( \{0\} \times B) = 1$, that is $$\nu_{B}( \{0\} \times A) = \frac{\nu(A \cap B)}{\nu(B)}.$$ Let $$\mu((0,s)) = \delta_{s}$$ for $s \in B$.  The decision problem is not conditionally rational, and the set of strategies that are not best responses to their assigned measures has $\nu_{B}$-measure one.
So each such problem has a profitable Borel deviation which gives rise to a partial uniformization of $R$, defined on a Borel set $\bar{B} \subseteq B$ of positive measure.  These uniformizations can be cobbled together into a Borel uniformization defined almost everywhere.

To see this, let
\[
\begin{aligned}
M=\sup\{\nu(C)\mid {}&
C \text{ is Borel and}\\
& R \text{ has a Borel uniformization on } C\}.
\end{aligned}
\]
Using countable choice, choose pairs $(C_n,u_n)$ such that $\nu(C_n)$ converges to $M$. Set
\[
C=\bigcup_{n<\omega}C_n,
\]
and define $u$ on $C$ by using $u_n$ on
$C_n\setminus\bigcup_{m<n}C_m$.
Then $u$ is a Borel uniformization of $R$ on $C$.
By the definition of $M$ and the choice of the $C_n$, we must have $\nu(C)=M$.

If $M<1$, the measurable set $2^{\omega} \setminus C$
contains a Borel set $B$ of positive measure.
The preceding argument gives a Borel set $E\subseteq B$
of positive measure and a Borel uniformization $v$ on $E$.
Then $u\cup v$ uniformizes $R$ on $C\cup E$, whose
measure exceeds $M$, a contradiction.

Hence $\nu(C)=1$.
Extending $u$ by a fixed value outside $C$ gives a Borel
map $f:2^\omega\to 2^\omega$ satisfying
\[
\nu\bigl(\{x\in C \ | \ (x,f(x))\notin R\}\bigr)=0.
\]
This proves $\mathsf{MUP}^{\mathsf{UM}}$. \end{proof}

\begin{prop}  Assume $\mathsf{V} = \mathsf{L}$. There is an analytically measurable decision problem which is deviation stable but not conditionally rational.
\end{prop}
\begin{proof}  This follows from Proposition 3.9 and Proposition 2.4(7). Let $R$ be the graph of a total coanalytic function which does not contain a perfect subset.  We can use the decision system as in Prop. 3.9 with $B = 2^{\omega}$.  The map $\pi^R$ is analytically measurable, and any Borel partial uniformization yields a contradiction as uncountable Borel sets contain perfect subsets.  \end{proof}

\begin{prop} There is a universally measurable decision problem which is deviation stable but not conditionally rational.
\end{prop}
\begin{proof} We will use the same template from Definition 3.8 with the graph of a function with nonmeasurable domain and universally null range as in Proposition 2.4 (5). Let $\lambda$ be the coin-flip measure on $2^{\omega}$.  Using a measure preserving Borel isomorphism to move to $2^{\omega}$, we can apply a theorem of Grzegorek (see Theorem 2.8 of \cite{LNSh}) to find a universally null set of the same cardinality as a nonmeasurable set. Thus we may choose a nonmeasurable set
$A\subseteq 2^{\omega}$, a universally null set
$N\subseteq 2^{\omega}$, and a bijection $f:A\to N$.
Write
\[
F=\{(s,f(s)) \ | \ s\in A\}.
\]
As in Proposition 2.4(5), $F$ is universally null and hence universally measurable.  Let $\pi = \pi^F$ as in Definition 3.8.  Thus $X = \{0,1\} \times 2^{\omega}$ and $Y = 2^{\omega}$.  The payoff $\pi$ is universally measurable. Let
\[
P=(s\mapsto(0,s))_*\lambda,
\qquad
\mu((i,s))=\delta_s.
\]
The assignment $\mu$ is continuous, and 
$$\pi( (0,s), \mu( (0,s)) ) = 0$$ for a set of $(0,s)$ of $P$-measure one.  For any Borel deviation $$d:X\to X,$$ write
\[
d((0,s))=(i(s),g(s)).
\]
The key point is $$\lambda\{ s \in A \ | \ g(s) = f(s)\} = ((\text{id},g)_*\lambda)(F) = 0$$ by the fact that $F$ is universally null and $$(\text{id},g)_*\lambda \in \Delta(2^{\omega} \times 2^{\omega})$$ is an atomless Borel probability measure as $g$ is Borel. Therefore, setting $E = \{s in A \ | \ g(s) = f(s) \}$,
\[
\int_X\pi(d(x),\mu(x))\,dP(x)
= 
\]
\[
\int_E \pi((i(s),g(s)),s) d\lambda(s) +  \int_{2^{\omega} \setminus E}  \pi((i(s),g(s)),s) d\lambda(s) = 0.
\]
Thus there are no profitable Borel deviations.  On the other hand, for every $s\in A$,
\[
\pi((1,f(s)),\delta_s)=1
>
0=\pi((0,s),\delta_s),
\]
and for $s\notin A$, every strategy earns zero against
$\delta_s$. Hence the set of nonrational strategies
on the support of $P$ is exactly $\{0\}\times A$.
Since $A$ is not measurable it is not contained in any Borel null set, so the problem is not conditional rational.
\end{proof}

We can now address the two issues (P1) and (P2)  simultaneously by requiring of our decision systems that the set of strategies at which rationality holds is measurable by the prior.  This rules out the pathology of Proposition 3.11, so the set we want to uniformize is universally measurable with measurable projection. 
 
\begin{defn}  We say a decision problem $(\pi,\mu,P)$ is rationality-measurable if $$\{x \ | \ x \mbox{ is not a best response to }\mu(x) \}$$ is $P$-measurable. 
\end{defn}

\begin{rem}  Note that $\mu$ is a Borel map, so we are requiring that $$\{ x \ | \ (x,\mu(x)) \in E^{\pi} \}$$ be $P$-measurable where $$E^{\pi} = \{ (s,\nu) \ | \ s \mbox{ is a best response to }\nu \}$$ is the best response relation for $\pi$.  So rationality-measurable is weaker than requiring that $E^{\pi}$ be universally measurable. \end{rem}

\begin{thm}  The following are equivalent:
\begin{enumerate}
    \item Every rationality-measurable decision problem which is deviation stable is conditionally rational.    
    \item $\mathsf{MUP}^{\mathsf{UM}}$
\end{enumerate}
\end{thm}

\begin{proof}  The rationality-measurable condition implies that $$B= \{ (s,a) \ | \ \pi(a,\mu(s)) > \pi(s,\mu(s))\},$$ which is universally measurable, has $P$-measurable projection so we can find a Borel $C \subseteq \mbox{proj}(B)$ with $P(C) > 0$ and apply 
$\mathsf{MUP}^{\mathsf{UM}}$ to $\bar{B} = B \cap (C \times X)$.  For the other direction, the decision problems we considered in Proposition 3.9 had a $P$-measure one failure set as the relation was total.  Thus that proof works under the present restriction. \end{proof}


\section{Games}

We study common priors in the context of games, and the results from the last section will carry over by considering the decision problems $$(\pi_i,\mu^i,P_i)$$ where $\pi_i$ is the payoff of player $i$, $P_i$ is the marginal of the common prior $P$, and $\mu^i$ is a disintegration of $P$.

\begin{defn}[Common Strategy Priors and Disintegration]
Suppose $X_1$ and $X_2$ are Polish spaces.  We call a probability measure
$P\in\Delta(X_1\times X_2)$ a \emph{common strategy prior}, and write $$P_i = \mbox{marg}_{X_i}(P)$$ for its marginal on $X_i$.
For each $i$, there is a Borel map
\[
\mu^i:X_i\longrightarrow\Delta(X_j)
\]
such that, for all Borel sets $D_1\subseteq X_1$ and
$D_2\subseteq X_2$,

\begin{equation}
P(D_1\times D_2)
=
\int_{D_i}\mu^i(s)(D_j)\,dP_i(s).
\end{equation}\label{eq:generatesprior} 
Moreover, $\mu^i$ is unique up to $P_i$-almost-everywhere
equality (see Theorems 5.3, 5.4 of \cite{kallenberg}).  If there is an underlying game we say $P$ is conditionally rational if $x$ is a best response to $\mu^i(x)$ on a set of $P_i$-measure one for each $i$.
\end{defn}

\begin{defn}[Self-supported Measure Assignment]  
Suppose $G$ is a bounded, universally measurable game.
We say that a pair of partial maps $(\mu^1,\mu^2)$ is a rational self-supported measure assignment if $$(s,\mu^i(s)) \in E^G_i \mbox{ and } \mu^i(s)(A_j) = 1$$ for every $s \in A_i$, where $A_j = \mbox{dom}(\mu^j)$ and $E^G_i$ is the best response relation for player $i$. \end{defn}

\begin{prop}[Conditionally Rational]  Suppose $G$ is a bounded, universally measurable game and $P \in \Delta(X_1 \times X_2)$ is a common strategy prior.  Then the following are equivalent:

\begin{enumerate}
  \item  $P$ is conditionally rational
  \item There are Borel sets $A_i \subseteq X_i$ with $P_i(A_i) = 1$ and a rational self-supported measure assignment on $A = A_1 \times A_2$ which generates 
  $P$ in that $$P(B_1\times B_2) = \int_{A_1 \cap B_1}\mu^1(s)(B_2)\,dP_1(s) = \int_{A_2 \cap B_2}\mu^2(s)(B_1)\,dP_2(s)$$ for every Borel $B_i \subseteq X_i$.

\end{enumerate}
\end{prop}
\begin{proof} For the forward direction, let $A^0_i$ be a $P_i$-full measure Borel set on which every $s$ is a best reponse to $\mu^i(s)$ and recursively define Borel sets
\[
A_i^{n+1}
=\{s\in A_i^n:\mu^i(s)(A_j^n)=1\}.
\]
The conditional probability identities imply inductively
that each $A_i^n$ is Borel and has full marginal measure.
Consequently, the Borel sets
\[
A_i=\bigcap_{n<\omega}A_i^n
\]
satisfy $P(A_1\times A_2)=1$ and, for every $s\in A_i$,
\[
(s,\mu^i(s))\in E_i^G,
\qquad
\mu^i(s)(A_j)=1.
\]
For the other direction we can extend each $\mu^i$ using an arbitrary constant value outside of $A_i$.
\end{proof}

\begin{defn}[Correlated Equilibrium]  A conditionally rational common strategy prior on a game is called a correlated equilibrium of the game.\end{defn}

\begin{rem} Thus correlated equilibrium can be characterized by a rational self-supported measure assignment on a Borel rectangle together with a pair of marginal probability measures that give rise to a common prior through the formula of equation (4.1).  We note that the existence of such measure assignments implies the strategy rectangle has the best response property in the sense of Pearce \cite{pearce}, and thus correlated equilibria concentrate on rationalizable strategies.
\end{rem}

\begin{defn}[Deviation Stable Common Strategy Prior] Suppose $G$ is a game with bounded universally measurable payoffs and Polish strategy spaces $X_i$, and $P \in \Delta(X_1 \times X_2)$ is a common prior.  We say $P$ is a deviation stable common strategy prior of $G$ if for each Borel map $d:X_1 \rightarrow X_1$, $$\int \pi_1(s,t) dP \geq \int \pi_1(d(s),t) dP,$$ and similarly for player 2.\end{defn}

\begin{rem}  This formulation has its origin in Aumann \cite{aumann1} where it is used as the definition of  correlated equilibrium.  See also Definition 2.1 of \cite{aumann2} and Definition 2.7 of \cite{stein}.
It has the virtue of only referencing the common strategy prior.  Also, it is easy to see that if $$P = P_1 \otimes P_2 \in \Delta(X_1 \times X_2)$$ is a common prior, then $P$ is deviation stable if and only if $(P_1,P_2)$ is a Nash equilibrium.  The conditionally rational formulation is derived as a consequence in the finite setting (see for example Proposition 2.3 of \cite{aumann2}).  On the other hand, in Yildiz \cite{Yildiz}, \emph{Graduate Game Theory}, Definition 2.13 and Theorem 2.7, the conditionally rational formulation is used as the definition of correlated equilibrium, which is the approach we take here.
\end{rem}

\begin{prop} For a game $G$ with Polish strategy spaces and bounded universally measurable payoffs, if a common strategy prior $P$ is a correlated equilibrium of $G$, then it is deviation stable. \end{prop}

\begin{proof}  Immediate.
\end{proof}

\begin{prop}  Assume $\mathsf{MUP}$.  For a game $G$ with Polish strategy spaces and bounded universally measurable payoffs, and a common strategy prior $P$, the following are equivalent:

\begin{enumerate}
    \item $P$ is a correlated equilibrium of $G$
    \item $P$ is a deviation stable.
\end{enumerate}
\end{prop}

\begin{proof} Let $\mu^i$ be a disintegration of $P$, and $P_i$ be the marginals of $P$.  Then we may apply Proposition 3.7(3) to the decision problems $(\pi_i,\mu^i,P_i)$, where $\mu^i$ is a Borel disintegration of $P$ and $P_i$ is its marginal for player $i$.  The other direction is immediate.
\end{proof}

\begin{prop}[Uniformization Necessary]
Assume that every deviation stable common strategy prior is conditionally
rational in every game with compact Polish strategy spaces and bounded universally measurable payoffs.
Then $\mathsf{MUP}^{\mathsf{UM}}$ holds.
\end{prop}

\begin{proof}  This follows from the decision theoretic result in Proposition 3.9. We can can set player 2's payoff equal to zero, and use the prior $$Q(A \times B) = \int_A \mu(x)(B) dP(x)$$ for Borel sets $A \subset X$ and $B \subseteq Y$, where $(\pi,\mu,P)$ is one of the relevant decision problem. \end{proof}

\begin{rem}  As a corollary, by Proposition 2.4(6), if every deviation stable common strategy prior is conditionally rational for universally measurable games then $\bDelta^1_2 \subseteq \Gamma^{\mathsf{UM}}$.  This latter hypothesis is shown by Levy in \cite[Theorem~3]{Levy} to be equivalent to the existence of measurable selectors for Nash equilibrium of parameterized Borel games.  We observe here that his conclusion follows from $\mathsf{MUP}_\mathsf{UM}$ directly.  $\mathsf{MUP}_\mathsf{UM}$ implies that any total function with universally measurable graph is universally measurable.  The equilibrium relation in the setting of \cite{Levy} is coanalytic and has full projection so its Kondo uniformization must be universally measurable.  \end{rem}

\begin{thm}  The principle $\mathsf{MUP}^{\mathsf{UM}}$ is equivalent to the assertion that any deviation stable rationality-measurable common strategy prior is a correlated equilibrium in any game with bounded, universally measurable payoffs and Polish strategy spaces. \end{thm}
\begin{proof} By rationality measurable we mean that each marginal $P_i$ of the common prior measures the set of strategies which are best responses to their assigned measures (see Definition 3.12).  The result follows from Theorem 3.14. \end{proof}

\begin{thm} Deviation stable priors are conditionally rational for every game with bounded, Borel measurable payoffs and Polish strategy spaces, while the analogous statement for analytically measurable games is independent of $\mathsf{ZFC}$.  \end{thm}
\begin{proof}  The Borel case follows from Prop. 3.7.
We can use $\mathsf{L}$ for failure in the analytically measurable case by Prop. 3.10.  We note that the decision problem for Proposition 3.10 associated to the graph $R$ of the coanalytic map from Proposition 2.4(7) is analytically measurable, as the preimage of any open set is either analytic or coanalytic.  For the other part we can use a forcing extension of $\mathsf{L}$ in which $\mathsf{MA} + \neg \mathsf{CH}$ holds by Prop. 2.4 (2).  \end{proof}

\section{Epistemics}

We use the universal type structure from Brandenburger and Dekel \cite{brand}.  For Polish spaces $X_i$ there are Polish spaces $T_i$, called the universal type spaces, and homeomorphisms $$\lambda_i:T_i \rightarrow \Delta(X_j \times T_j).$$  If $G = (\pi_1,\pi_2)$ is a game on $X_1,X_2$ we say a strategy $s \in X_i$ is rational at $t \in T_i$ if $s$ is a best response to $\mbox{marg}_{X_j}(\lambda_i(t))$ and let $R_i^1(G)$ denote the set of such pairs for each $i$.  Define for $n\geq 1$,
$$R_i^{n+1}(G) = \{ (s,t) \in R_i^n(G) \ | \ \lambda_i(t)(R_j^n(G)) = 1\}, $$ and $$\text{RCBR}_i(G) = \mbox{proj}_{X_i}( \bigcap_{n > 0} R_i^n(G)).$$
\noindent We set $$\text{RCBR}(G) =  \text{RCBR}_1(G) \times \text{RCBR}_2(G),$$ and we say $s \in \text{RCBR}_i(G)$ is consistent with rationality and common belief in rationality.  Here we have used the convention that when $i$ is used to index a player, $j$ will represent the other player.

There is a subtle point here.  Even if the best response relation $E_i^G \in \Gamma^{\mathsf{UM}}$, we have  $\text{RCBR}_i(G) \in \exists^{\mathbb{R}} \Gamma^{\mathsf{UM}}$.  But as we identify measures with their completions, the sequence defined above is well-defined, as is the statement of the theorem.

\begin{lem}[Type Lift Lemma]  Suppose $A_i \in \Gamma^{\mathsf{UM}}$ for each $i$ and that the universally measurable maps $$\mu^i:A_i \rightarrow \Delta(X_j)$$ form a rational self-supported measure assignment on $A_1 \times A_2$.  Then there are universally measurable type assignments $$t^i:A_i \rightarrow T_i$$ such that for each $i$ and $s \in A_i$, $$\mbox{marg}_{X_j}(\lambda_i(t^i(s))) = \mu^i(s) \mbox { and } \lambda_i(t^i(s))(\mbox{graph}(t^j))= 1,$$ equivalently,\footnote{Note that a universally measurable function function between Polish spaces has a universally measurable graph.}
$$\lambda_i(t^i(s)) = (\mbox{id}_{A_j},t^j)_\ast \mu^i(s).$$ \end{lem}

\begin{proof}  The lift argument involves the construction of the type spaces $T_i$.  To sketch it, set $X_i^0 = X_j$, $t^i_1(s) = \mu^i(s)$ for each $i$ and 
inductively for $n \geq 0$ define 
$$X_i^{n+1} = X_i^{n} \times \Delta(X_j^{n}),$$
and for $n \geq 1$ define $$t_{n+1}^i(s) = (\mbox{id}_{A_j},t_{1}^j,...,t_{n}^j)_\ast \mu^i(s) \in \Delta(X_i^n)$$
simultaneously for each $i$.  Each $t_{n}^i$ is a universally measurable function by induction, since pushing forward a universally measurable measure assignment by a universally measurable map is universally measurable. Thus $$t^i(s) = (t_{1}^i(s),t_{2}^i(s), t_{3}^i(s),...) \in T_i,$$ which is easily seen by induction to satisfy common certainty of coherency, is the desired type assignment. By uniqueness of the Kolmogorov extension $$\lambda_i(t^i(s)) =(\mbox{id}_{A_j},t^j)_*\mu^i(s)$$ as they have the same finite dimensional marginals (see \cite{brand} and \cite{Zoble} for more details).
\end{proof}

\begin{thm}[Correlated Equilibrium is RCBR] Suppose $G$ is a game and $P \in \Delta(X_1 \times X_2)$ is a correlated equilibrium. Then $$P(\textnormal{RCBR(G)}) = 1.$$ \end{thm}

\begin{proof}  That is, a conditionally rational equilibrium concentrates on strategies which are consistent with rationality and common belief in rationality.  Let $P \in \Delta(X_1 \times X_2)$ be such an equilibrium.   Then by Proposition 4.3 there is a disintegration of $P$ which is a rational self-supported Borel measure assignment $(\mu^1,\mu^2)$ on some Borel strategy rectangle $A_1 \times A_2$ with 
$$P(A_1 \times A_2) = 1.$$  
By Lemma 5.1, the type lifts $$t^i:A_i \rightarrow T_i,$$
satisfy
\[
\lambda_i(t^i(s))
=(\operatorname{id}_{A_j},t^j)_*\mu^i(s) \in \Delta(X_j \times T_j),
\]
where the measure $\mu^i(s)$ is regarded as a probability
measure on $A_j$.  We note that as each $\mu^i$ is a Borel function, so is each $t^i$.

Thus for $s \in A_i$, $$\lambda_i(t^i(s))( \mbox{graph}(t^j)) = 1$$ and $$\mbox{marg}_{X_j}(\lambda_i(t^i(s))) = \mu^i(s).$$
It follows immediately by induction that for $n \geq 1$, $$(s,t^i(s)) \in R_i^n$$ for each $i$ and $s \in A_i$, so these type assignments witness
\[
A_1\times A_2\subseteq \text{RCBR(G)}, 
\]
and thus $P(\text{RCBR(G)}) = 1$.  \end{proof}

\begin{thm}
Assume $\mathsf{V}=\mathsf{L}$. There is a bounded, analytically measurable game $G$
on Polish strategy spaces with a deviation stable common strategy prior $P$ such that $$P(\mathrm{RCBR}(G)) = 0.$$ 
\end{thm}

\begin{proof} We use the argument of Proposition 3.10 with additional strategies added for player 1.  Let $$R \subset 2^{\omega} \times 2^{\omega}$$ be the graph of a total function which is coanalytic and does not contain a perfect subset, and let $\lambda \in \Delta( 2^{\omega})$ denote the coin-flip measure. Thus every section $R_x$ is nonempty,
while $$\lambda (\{x \ | \ (x,g(x))\in R\})=0$$
for every Borel map $g:2^{\omega} \to 2^{\omega}$.  Set $X_1= \{0,1,2\} \times 2^{\omega}$ and $X_2=2^{\omega}$, and define
\[
\pi_1((n,a),x)=
\begin{cases}
  - 1  & n=0 \wedge x \neq a,\\
   0  & n=0 \wedge x = a,\\
  \chi_R(x,a) & n=1,\\
  - 1/2 & n = 2 \wedge x \neq a, \\
   0  & n = 2 \wedge x = a. 
\end{cases}
\]
We give player $2$ the constant payoff $\pi_2=0$ and form the game $G = (\pi_1,\pi_2)$.
These payoffs are bounded and analytically measurable.
Let $P$ be the pushforward of $\lambda$ under
\[
  x\longmapsto ((0,x),x).
\]
Both players receive payoff zero under $P$.
For any Borel deviation of player $1$, write its restriction
to part $0$ as
\[
  d_1(0,x)=(n(x),a(x)).
\]
The functions $n$ and $a$ are Borel. The deviating payoff is
nonpositive when $n(x)\ne1$, and it can be positive when
$n(x)=1$ only if $(x,a(x))\in R$, which occurs on a
$\lambda$-null set. Hence no Borel deviation is profitable.
Player $2$ is indifferent, so $P$ is deviation stable. 

To show that $P$ concentrates on strategies that are not rationalizable (and hence not in $\text{RCBR}(G)$), fix $(0,a)\in X_1$ and $\nu\in\Delta(2^{\omega})$.
If $\nu\ne\delta_a$, then
\[
  \pi_1((2,a),\nu)
  =-\frac12\bigl(1-\nu(\{a\})\bigr)
  >-\bigl(1-\nu(\{a\})\bigr)
  =\pi_1((0,a),\nu).
\]
If $\nu=\delta_a$, choose $b\in R_a$. Then
\[
  \pi_1((1,b),\delta_a)=1
  >0=\pi_1((0,a),\delta_a).
\]
Thus no strategy in $\{0\}\times 2^{\omega}$ is a best response to a belief. Consequently, the Borel rectangle
\[
  A=(\{0\}\times 2^{\omega})\times 2^{\omega}
\]
satisfies $P(A)=1$ and $A\cap\mathrm{RCBR}(G)=\varnothing$. \end{proof}

Theorem 4 from the introduction now follows.  Under the hypothesis $\mathsf{MA} + \neg \mathsf{CH}$, every deviation stable common strategy prior for an analytically measurable game is conditionally rational by Theorem 4.13 and hence concentrates on $\mathrm{RCBR}$ strategies by Theorem 5.2.  Theorem 5.3 gives a counterexample in $\mathsf{L}$.

\begin{defn}[Common Prior on States] Let
\[
\Omega=(X_1\times T_1)\times(X_2\times T_2),
\]
where $\lambda_i:T_i\to\Delta(X_j\times T_j)$ are the
universal type maps.  We call $\Omega$ the set of (joint) states.  For a state $\omega = ((s_1,t_1),(s_2,t_2)) \in \Omega$ we can call $(s_i,t_i)$ an individual state or a state for player $i$. For a Borel event $B\subseteq\Omega$ 
write $B_{s_i,t_i}$ for its section in the opponent's
coordinates. Thus player $i$'s probability of $B$ at
a state $\omega$ is
\[
\lambda_i(t_i)(B_{s_i,t_i}).
\]
A probability measure $P\in\Delta(\Omega)$ is a
\emph{common prior} if, for each player $i$, the
conditional distribution of the opponent's strategy
and type given $(s_i,t_i)$ is $\lambda_i(t_i)$,
$P$-almost surely. Equivalently, writing $P_i$ for
the marginal on $X_i\times T_i$,
\[
P(B)=
\int_{X_i\times T_i}
\lambda_i(t_i)(B_{s_i,t_i})\,P_i(d s_i,d t_i)
\]
for every Borel $B\subseteq\Omega$.
\end{defn}

\begin{defn}  An \emph{epistemic correlated equilibrium} is a common prior $Q$ which concentrates on states at which rationality and common belief in rationality hold.  \end{defn}

\begin{prop}  Every correlated equilibrium $P$ lifts to an epistemic correlated equilibrium $P^*$ with $$\mbox{marg}_{X_1 \times X_2}(P^*) = P.$$ Conversely, the strategy marginal $\mbox{marg}_{X_1 \times X_2}(P^*)$ of an epistemic correlated equilibrium is a correlated equilibrium of the underlying game. \end{prop}

\begin{proof}  Let $P$ be a correlated equilibrium.  The lift is defined by
$$P^*=F_*P,$$
where $$F: A_1 \times A_2 \rightarrow \Omega$$ is given by 
$$F(s_1,s_2) = (F_1(s_1),F_2(s_2))
=\bigl((s_1,t^1(s_1)),(s_2,t^2(s_2))\bigr),$$
and $t^i(s)$ is the type lift from Theorem 5.1.
Thus every state in the image of $F$ satisfies rationality and common belief in rationality.
\noindent Each set 
\[
C_i=F_i(A_i)
=
\{(s,t^i(s)) \ | \ s\in A_i\}
\]
is Borel and the lifting identity implies that
\[
\lambda_i(t^i(s))(C_j)=1
\]
and
\[
\operatorname{marg}_{X_j}\lambda_i(t^i(s))
=
\mu^i(s).
\]
Thus every point of $C_i$ is rational and assigns
probability one to $C_j$. Consequently $C_1\times C_2$
is self-supporting and $$P^*(C_1\times C_2)=1.$$

A disintegration of $P^*$ over its player $i$ marginal $P_i^*$
is given by $\lambda_i(t_i)$ for $P_i^*$-almost every
$(s_i,t_i)$. This follows from the disintegration of $P$
and the definition $P^*=F_*P$. Thus $P^*$ is a common prior on state space and hence an epistemic
correlated equilibrium.  

For the second assertion, assume $P^*$ is an epistemic
correlated equilibrium and $P$
its strategy marginal.  Intuitively, since $P^*$ is a common prior, each player's belief about their opponent's strategies, given their own strategy, is a mixture of marginals of associated types.  Since $P^*$ concentrates
on $\text{RCBR}$ states, each such strategy is a best response to almost every associated type's marginal, and is therefore a best response to the mixture. 

Let $P_i^* = \mbox{marg}_{X_i\times T_i}(P^*)$ and 
$P_i = \mbox{marg}_{X_i}(P^*_i)$.  For each $i$ let
$$\nu_i:X_i \rightarrow \Delta(T_i)$$ be a disintegration of $P_i^*$ over $P_i$. Also let 
\[
b_i(t_i)
=
\mbox{marg}_{X_j}(\lambda_i(t_i)).
\]

Since $P^*$ is a common prior, a disintegration of $P^*$
over $P_i^*$ is given by the type structure assignment $\lambda_i(t_i)$.  Thus, a disintegration $\mu^i$ of $P$ over $P_i$
is obtained by averaging the strategy beliefs as follows.
For each Borel set $B\subseteq X_j$, define

$$\mu^i(s_i)(B) = \int_{T_i} b_i(t_i)(B)\,d\nu_i(s_i).$$
This defines a probability measure on $X_j$, which we abbreviate as $$
\mu^i(s_i) = \int_{T_i} b_i(t_i)\,d\nu_i(s_i).$$
Because $P^*$ concentrates on RCBR states, it concentrates
on states at which each player is rational.  Hence, for
$P_i$-almost every $s_i$,
\[
s_i\in\operatorname{BR}_i
\bigl(\operatorname{marg}_{X_j}\lambda_i(t_i)\bigr)
\qquad
\text{for $\nu_i(s_i)$-almost every $t_i$}.
\]
Fix such an $s_i$. For every alternative $a_i\in X_i$,
integrating the best-response inequalities yields
\[
\begin{aligned}
\pi_i(s_i,\mu^i(s_i))
&=
\int_{T_i}
\pi_i\bigl(s_i,\operatorname{marg}_{X_j}\lambda_i(t_i)\bigr)
\,d\nu_i(s_i)\\
&\geq
\int_{T_i}
\pi_i\bigl(a_i,\operatorname{marg}_{X_j}\lambda_i(t_i)\bigr)
\,d\nu_i(s_i)\\
&=
\pi_i(a_i,\mu^i(s_i)).
\end{aligned}
\]
Thus $s_i\in\operatorname{BR}_i(\mu^i(s_i))$ for
$P_i$-almost every $s_i$. Since this holds for each player,
$P$ is conditionally rational. \end{proof}

As a concluding remark, we noted in Remark 4.5 that a rational self-supported measure assignment on a rectangle $A_1 \times A_2$ is a witness for the rectangle having the best response property with respect to a game $G$.  The maximal strategy rectangle of $G$ with the best response property exists, and constitutes the rationalizable strategies $$\text{RAT}(G) = \text{RAT}_1(G) \times \text{RAT}_2(G).$$  Thus $$\text{RCBR}(G) \subseteq \text{RAT}(G)$$ as $\text{RCBR}(G)$ has the best response property.  Assuming $\mathsf{AD}_{\mathbb{R}}$, the Axiom of Real Determinacy, for an arbitrary game $G$ played on Polish strategy spaces we can find a uniformization $\mu^i$ of the best response relation defined on $\text{RAT}_i(G)$ for each $i$, and consisting only of measures that concentrate on $\text{RAT}_j(G)$.  As all subsets of Polish spaces are universally measurable under this hypothesis, the pair 
$(\mu^1, \mu^2)$ is a universally measurable rational self-supported measure assignment on $\text{RAT}(G)$.  Hence Lemma 5.1 applies (together with the induction of Theorem 5.2) and we conclude that $$\text{RCBR}(G) = \text{RAT}(G)$$ for this arbitrary game, yielding the maximal extension of the Fundamental Theorem of Epistemic Game Theory under determinacy.  These ideas are developed in \cite{Zoble}.

\end{document}